\documentclass{amsproc}
\usepackage{xypic}
\usepackage{amssymb}
\usepackage{amsmath}

\newtheorem{theorem}{Theorem}[section]

\theoremstyle{definition}
\newtheorem{definition}[theorem]{Definition}
\newtheorem{conjecture}[theorem]{Conjecture}
\newtheorem{example}[theorem]{Example}
\newtheorem{examples}[theorem]{Examples}
\newtheorem{proposition}[theorem]{Proposition}

\newtheorem{corollary}[theorem]{Corollary}
\renewcommand{\to}{\xymatrix@1@=15pt{\ar[r]&}}
\renewcommand{\mapsto}{\xymatrix@1@=15pt{\ar@{|->}[r]&}}

\newcommand{\dual}{^\vee}
\newcommand{\km}{{\mathcal M}}
\newcommand{\ko}{{\mathcal O}}
\newcommand{\kp}{{\mathcal P}}
\newcommand{\kc}{{\mathcal C}}
\newcommand{\kk}{{\mathcal K}}
\newcommand{\ke}{{\mathcal E}}
\newcommand{\kf}{{\mathcal F}}
\newcommand{\kl}{{\mathcal L}}
\newcommand{\kt}{{\mathcal T}}

\newcommand{\IC}{{\mathbb C}}

\newcommand{\IP}{{\mathbb P}}
\newcommand{\IQ}{{\mathbb Q}}
\newcommand{\IG}{{\mathbb G}}
\newcommand{\IR}{{\mathbb R}}
\newcommand{\IZ}{{\mathbb Z}}

\newcommand{\Pic}{{\rm Pic}}
\newcommand{\Hom}{{\rm Hom}}
\newcommand{\Br}{{\rm Br}}
\newcommand{\Ext}{{\rm Ext}}

\newcommand{\Coh}{{\rm\bf Coh}}
\newcommand{\Db}{{\rm D}^{\rm b}}

\newcommand{\cal}{\mathcal}
\newcommand{\ka}{\cal A}
\def\mydate{\number\day\space\ifcase\month \or January\or February\or March\or
April\or May\or June\or July\or August\or September\or October\or
November\or December\fi \space\number\year}

\theoremstyle{remark}
\newtheorem{remark}[theorem]{Remark}

\numberwithin{equation}{section}

\begin{document}
\title[The Global Torelli Theorem: classical, derived, twisted.]
      {The Global Torelli Theorem: classical, derived, twisted.}

\author[Daniel Huybrechts]{Daniel Huybrechts}

\address{Mathematisches Institut, Universit\"at Bonn,
Beringstr.\ 1, D - 53115 Bonn, Germany}

\email{huybrech@math.uni-bonn.de}




\maketitle

These notes survey work on various aspects and generalizations of
the Global Torelli Theorem for K3 surfaces done over the last ten
years. The classical Global Torelli Theorem was proved a
long time ago (see \cite{BR,LP,PS,Siu}), but the interest in
similar questions has been revived by the new approach to K3
surfaces suggested by mirror symmetry.

Kontsevich proposed  to view mirror symmetry as an equivalence
between the bounded derived category of coherent sheaves on a
Calabi--Yau manifold and the derived Fukaya category of its mirror dual.
For an algebraic geometer the bounded derived category of coherent
sheaves on a variety is a familiar object, but to view it as an
interesting invariant of the variety rather than a technical tool
to deal with cohomology is rather surprising. Due to results of
Mukai, Orlov, and Polishchuk, the bounded derived category of coherent sheaves
on an abelian variety is completely understood, i.e.\ we know
when two abelian varieties yield equivalent derived categories
and what the group of autoequivalences looks like.

Independently of their importance in mirror symmetry, K3 surfaces
form the next most simple class of Calabi--Yau manifolds and one
would  like to study them from a derived category point of view,
too. The program has been started already by Mukai in \cite{Mu}
and the break-through came with Orlov's article \cite{Or}. But
this was still not the end. For many reasons (mirror symmetry
considerations, existence of non-fine moduli spaces, geometric
interpretation of conformal field theories, etc.) one would like
to incorporate B-fields or, closely related, Brauer classes in the
picture. These notes will mostly concentrate on aspects that are
related to generalizations of the Global Torelli Theorem in this
direction. The following list of topics gives an idea of what
shall be discussed:
\smallskip

$\bullet$ Hitchin's generalized Calabi--Yau structures.

$\bullet$ The period of a twisted K3 surface.

$\bullet$ Brauer group and B-fields.

$\bullet$ Derived categories of twisted sheaves.

\bigskip


This survey contains essentially no proofs. I  have tried to
emphasize ideas and  refer to the literature for details. Some of
the results can very naturally be viewed in terms of moduli
spaces of generalized K3 surfaces and the action of the mirror
symmetry group, but I have decided to neglect these aspects almost
completely. Although mirror symmetry has shaped our way of
thinking about derived categories, the symplectic side of the
theory will not be touched upon.

\medskip

{\bf Acknowledgements:} Most of my own results on derived
equivalences surveyed here have been obtained jointly with Paolo
Stellari. I wish to thank him for the fruitful collaboration. I
am grateful to Jochen Heinloth, Max Lieblich, and Balazs
Szendr\H{o}i for comments on an earlier version.


\section{The classical Torelli}

A poetic explanation of the name `K3' was given in Andr\'e Weil's
remarks on his only article on the subject (which was in fact a
report for a U.S. Air Force grant). He writes:

\smallskip

{\it... il s'agit des vari\'et\'es k\"ahl\'eriennes dites K3,
ainsi nomm\'ees en l'honneur de Kummer, K\"ahler, Kodaira, et de
la belle montagne K2 au Cachemire.}

\medskip

The official definition goes as follows:

\begin{definition}
A connected compact complex surface $X$ is a \emph{K3 surface} if
its canonical bundle is trivial, i.e.\ $\omega_X\cong\ko_X$, and
$H^1(X,\ko_X)=0$.
\end{definition}

A trivializing section of $\omega_X$, i.e.\ a non-trivial
holomorphic two-form, will usually be denoted $\sigma\in
H^0(X,\omega_X)$. It is unique up to scaling.

\begin{examples}\label{ExK3}
 i) The \emph{Fermat quartic} is a concrete example of a (projective) K3 surface.
 It is given as the hypersurface in $\IP^3$ described by the
 equation $x_0^4+x_1^4+x_2^4+x_3^4=0$. The adjunction formula shows that the
 canonical bundle is trivial and standard vanishing results
 for the cohomology of line bundles on the projective space prove
 the required vanishing.

 ii) \emph{Kummer surfaces} form another distinguished type
 of K3 surfaces. If $T$ is an abelian surface or a complex torus
 of dimension two, then the associated Kummer surface ${\rm Kum}(T)$
 is the minimal resolution of the quotient $T/\pm$ by the standard involution
 $x\mapsto -x$ (which has sixteen fixed-points). In particular, ${\rm Kum}(T)$ contains sixteen
 $(-2)$-curves, i.e.\ smooth irreducible rational curves $C_i$
 which by adjunction satisfy $(C_i.C_i)=-2$. Note that ${\rm
 Kum}(T)$ is projective if and only if $T$ is projective.
\end{examples}

In the following we briefly recall a few standard facts from the
theory of K3 surfaces, for further details see \cite{Per,BPV}:

\smallskip

{\bf 1.} {\bf Any K3 surface is diffeomorphic to the Fermat
quartic.} This is shown by a simple deformation argument. As it
turns out, quartic K3 surfaces as well as Kummer surfaces are
arbitrarily close to any K3 surface. In particular, K3 surfaces
are simply-connected and the second cohomology $H^2(X,\IZ)$ is
therefore torsion free.

\smallskip

{\bf 2.} {\bf The intersection pairing} endows the middle
cohomology $H^2(X,\IZ)$ with the structure of a unimodular
lattice of rank $22$ which is abstractly isomorphic to
$$-E_8\oplus -E_8\oplus U\oplus U\oplus U.$$
Here, $E_8$ is the unique unimodular, positive-definite, even
lattice of rank eight and $U$ is the hyperbolic plane generated
by two isotropic vectors $e_1,e_2$ with $(e_1.e_2)=1$.

\smallskip

{\bf 3.} {\bf Any K3 surface is K\"ahler.} This is a deep theorem;
a complete proof of it was given by Siu in \cite{Siu}. The
analogous statement for higher dimensional (simply-connected)
holomorphic symplectic manifolds does not hold. Usually a
K\"ahler form will be denoted by $\omega$ and its K\"ahler class
by $[\omega]\in H^2(X,\IR)$. Although it will not be explicitly
mentioned anywhere in the text, Yau's result on the existence of
Ricci-flat metrics plays a central r\^ole in the theory. One way
to phrase it is to say that any K\"ahler class can be uniquely
represented by a K\"ahler form $\omega$ with
$\omega^2=\lambda(\sigma\bar\sigma)$ for some positive scalar
factor $\lambda$.

\smallskip

{\bf 4.} {\bf The weight-two Hodge structure}
$$H^2(X,\IC)=H^{2,0}(X)\oplus H^{1,1}(X)\oplus H^{0,2}(X)$$
on $H^2(X,\IZ)$ is of great importance in the study of K3
surfaces. By definition $H^{2,0}(X)\cong H^0(X,\omega_X)\cong
\IC\sigma$ is of dimension one. Moreover, $H^{0,2}(X)$ is complex
conjugate to $H^{2,0}(X)$ and $H^{1,1}(X)$ is orthogonal (with
respect to the intersection pairing) to $H^{2,0}(X)\oplus
H^{0,2}(X)$. Thus, the weight-two Hodge structure on the
intersection lattice $H^2(X,\IZ)$ of a K3 surface is determined by
the line $\IC\sigma\subset H^2(X,\IC)$.

The lattice $H^2(X,\IZ)$ together with the natural Hodge
structure of weight two is called the \emph{Hodge lattice} of $X$.
A \emph{Hodge isometry} between two lattices endowed  with
additional Hodge structures is by definition a lattice isomorphism
that respects both structures, the quadratic form of the lattices
and the Hodge structures.

In our geometric situation any integral $(1,1)$-class $\delta\in
H^2(X,\IZ)\cap H^{1,1}(X)$ with $(\delta.\delta)=-2$ induces a
Hodge isometry of the Hodge lattice $H^2(X,\IZ)$ given as the
reflection $s_\delta$ in the hyperplane $\delta^\perp$. More
explicitly, $$\xymatrix{s_\delta:\gamma\ar@{|->}[r]&
\gamma+(\gamma.\delta)\delta.}$$

\smallskip

{\bf 5.} {\bf The K\"ahler cone}  is the open cone
$$\kk_X\subset H^{1,1}(X,\IR):=H^2(X,\IR)\cap H^{1,1}(X)$$ formed by all
K\"ahler classes $[\omega]$.  As
$([\omega].[\omega])=\int_X\omega^2>0$, it is contained in one
connected component $\kc_X$ of the positive cone of all classes
$\gamma\in H^{1,1}(X,\IR)$ with $(\gamma.\gamma)>0$. (Since the
intersection pairing on $H^{1,1}(X,\IR)$ has signature $(1,19)$,
the only other connected component is $-\kc_X$.)

Conversely, a class $\gamma\in \kc_X$ is a K\"ahler class if and
only if $(\gamma.[C])>0$ for all $(-2)$-curves $C\subset X$. For
a higher-dimensional analogue see \cite{HuyKK}.

The hyperplanes $\delta^\perp$ orthogonal to integral
$(-2)$-classes $\delta\in H^{1,1}(X)$ cut $\kc_X$ in possibly
infinitely many chambers. If $\gamma\in \kc_X$ is a class in the
interior of one such chamber, then there exist $(-2)$-curves
$C_1,\ldots,C_n\subset X$ such that $s_{[C_1]}(\ldots
s_{[C_n]}(\gamma)\ldots)$ is a K\"ahler class. The K\"ahler cone
$\kk_X$ forms one chamber.

\bigskip

All these results are interwoven with the culmination of the
theory:

\begin{theorem}\label{CGT}{\bf (Classical Global Torelli)}
Two K3 surfaces $X$ and $X'$ are isomorphic if and only if there
exists a Hodge isometry $$H^2(X,\IZ)\cong H^2(X',\IZ).$$
\end{theorem}

The theorem has a long and complicated history. It has been proved
in various  degrees of generality by: Piateckii-Shapiro and
Shafarevich, Burns and Rapoport, Peters and Looijenga, and Siu.

\begin{remark}\label{RemGT}
i) Although the Global Torelli Theorem asserts that there is an
isomorphism $f:X\cong X'$ whenever there exists a Hodge isometry
$g:H^2(X,\IZ)\cong H^2(X',\IZ)$, a given Hodge isometry $g$ might
not be realized as a  Hodge isometry of the form $f_*$ for any
isomorphism $f$.

ii) One important example of a Hodge isometry $g$ that cannot be
realized as $f_*$ is provided by the reflection $s_{[C]}$
associated to a $(-2)$-curve $C\subset X$. Indeed, if $[\omega]$
is a K\"ahler class, then  $f_*[\omega]$ is also a K\"ahler class.
Hence, $(f_*[\omega].[C])=\int_Cf_*\omega>0$. On the other hand,
$(s_{[C]}[\omega].[C])=-([\omega].[C])<0$.

It can be shown that any given Hodge isometry $g:H^2(X,\IZ)\cong
H^2(X',\IZ)$ can be composed with finitely many reflections
$s_{[C_1]},\ldots, s_{[C_n]}$ associated to $(-2)$-curves
$C_i\subset X$, such that either the Hodge isometry
$g\circ\left(s_{[C_1]}\circ\ldots\circ s_{[C_n]}\right)$ or its
negative maps a K\"ahler class on $X$ to a K\"ahler class on $X'$.
This new Hodge isometry can then be lifted to a unique
isomorphism due to the following remark.

iii) The full Global Torelli Theorem  proves the following
assertion: If $g:H^2(X,\IZ)\cong H^2(X',\IZ)$ is a Hodge isometry
that sends at least one K\"ahler class on $X$ to a K\"ahler class
on $X'$, then $g$ is induced by a unique(!) isomorphism $f:X\cong
X'$.
\end{remark}

As has been mentioned before, any two K3 surfaces are
diffeomorphic. So, instead of viewing K3 surfaces as different
complex surfaces, they might be viewed as complex structures on a
specific differentiable manifold of (real) dimension four. (It is
known that any complex structure on $M$ does define the structure
of a K3 surface on $M$, see \cite{FM}.) To be more precise, we
denote by $M$ the differentiable manifold underlying the Fermat
quartic in i), Example \ref{ExK3}. (The complex structure of the
Fermat quartic induces a natural orientation of $M$, which we
will  fix throughout.)

Note, that in order to view a K3 surface $X$ as a complex
structure on $M$, one first needs to fix a diffeomorphism of the
differentiable manifold underlying  $X$ with our fixed manifold
$M$. The choice of the diffeomorphism is by no means unique, as
the diffeomorphism group of $M$ is far from
being trivial. Due to
Borcea \cite{Borcea} and Donaldson \cite{Donaldson} one knows
that the cohomology representation of the diffeomorphism group
yields a surjection
\begin{equation}\label{BorceaDonald}
\xymatrix{{\rm Diff}(M)\ar@{->>}[r]&{\rm O}_{+}(H^2(M,\IZ)).}
\end{equation}
Here ${\rm O}_+(H^2(M,\IZ))$
 is the group of all
lattice isomorphisms preserving the orientation of the positive
directions (see the explanations in Remark \ref{orientRem}).

Suppose now that a K3 surface $X$ is described by  a complex
structure $I$ on $M$. We write $X=(M,I)$. The holomorphic
two-form $\sigma$ (unique up to scaling) can be viewed as a
complex two-form $\sigma\in\ka_\IC^2(M)$ on $M$. This two-form
satisfies:$${\rm i)}~\sigma~ {\rm is~ closed, ~i.e.}~
~d\sigma=0,~~{\rm ii)}~\sigma\wedge\sigma\equiv0,{\rm~~~and~~~}
{\rm iii)}~\sigma\wedge\bar\sigma>0.$$

The last condition means that $\sigma\wedge\bar\sigma$ is a positive
multiple of the fixed orientation at every point of $M$.
The two-form $\sigma$ is also called the \emph{holomorphic volume form}
or the \emph{Calabi--Yau form} of $X=(M,I)$.

An easy observation, presumably due to Andreotti, shows that the
converse also holds. Indeed, any complex two-form
$\sigma\in\ka^2_\IC(M)$ satisfying i)-iii) is induced by a
complex structure in the above sense.  More explicitly, one
defines $T^{0,1}M$ as the kernel of $\sigma:T_\IC M\to T_\IC^* M$
and $T^{1,0}M$ as its complex conjugate. Conditions ii) and iii)
ensure that this results in a decomposition of $T_\IC M$ which
defines an almost complex structure. This almost complex
structure is integrable due to i).

The \emph{period} of a K3 surface $X=(M,I)$ is the point
$[\sigma]\in \IP(H^2(M,\IC))$, where $\sigma$ is the holomorphic
two-form that determines $I$.  The Local Torelli Theorem, stated
in \cite{W} and most likely due to Andreotti, roughly says that
any small change of the complex structure $I$ on $M$ up to
diffeomorphisms of $M$ is reflected by a change of the period,
i.e.\ of the line $\IC[\sigma]\subset H^2(M,\IC)$. This has led
Weil to conjecture a global version of this result and, in fact,
he gives two versions of it (see \cite{W}):

\begin{conjecture}
{\rm i)} Suppose $\IC[\sigma]=\IC [\sigma']\subset H^2(M,\IC)$.
Then there exists a diffeomorphism $f\in{\rm Diff}(M)$ isotopic
to the identity, i.e.\ $f$ is contained in the identity component
${\rm Diff_o}(M)$ of the full diffeomorphism group ${\rm Diff}(M)$,
such that $\sigma=f^*\sigma'$. \footnote{The expert reader of course knows
that this has no chance to be true as stated. One either
has to restrict to general complex structures $\sigma$ and
$\sigma'$ or to add the assumption that at least one K\"ahler
class with respect to $\sigma$ is also a K\"ahler class with
respect to $\sigma'$.}

{\rm ii)} If $ [\sigma]=g([\sigma'])$ for a lattice isomorphism
$g\in {\rm O}(H^2(M,\IZ))$, then there exists a diffeomorphism
$f\in {\rm Diff}(M)$ such that $\sigma=f^*\sigma'$.
\end{conjecture}

\begin{remark}
The second version has been established by the classical Global
Torelli Theorem as stated in its more algebro-geometric form in
Theorem \ref{CGT}, but i) is still open. So, the Global Torelli
Theorem for K3 surfaces has not been fully proven yet! In order
to deduce i) from ii) one would have to prove that ${\rm
Diff_o}(M)$ coincides with the kernel of the representation ${\rm
Diff}_*(M)$ of ${\rm Diff}(M)\to {\rm O}(H^2(M,\IZ))$.
\end{remark}

\begin{remark}\label{SurjPK3}
We conclude this section with the surjectivity of the period map.
The Global Torelli Theorem is equivalent to the assertion that the
\emph{period map}
$$\xymatrix{\kp:\{I\}/{\rm Diff}_*(M)\ar[r]&
Q,}~~\xymatrix{X=(M,I)\ar@{|->}[r]& [\sigma]}$$ is generically
injective. Here, $Q\subset \IP(H^2(M,\IC))$ is the \emph{period
domain}
$$Q=\{x\in \IP(H^2(M,\IC))~|~(x.x)=0,~(x.\bar x)>0\}.$$

Using the Global Torelli Theorem it has been proved
(see \cite{Per,Looij,Siu2,Tod}) that any $x$ in the period domain is the
period of some K3 surface $X=(M,I)$. In other words, the period
map $\kp$ is surjective. Although a Global Torelli Theorem does
not hold in higher dimensions, the surjectivity of the period map
could nevertheless be established in broader generality, see \cite{HuyHK}.
\end{remark}
\section{Generalized K3 surfaces}
In 2002 Hitchin \cite{Hit} introduced generalized complex and
generalized Calabi--Yau structures. Generalized Calabi--Yau
structures on K3 surfaces were in detail investigated in
\cite{HuyGen}.

If we think of a K3 surface as given by a Calabi--Yau form
$\sigma$ on the differentiable manifold $M$, then the following
definition is very natural

\begin{definition}
A \emph{generalized Calabi--Yau form} on $M$ is an even complex
form $\varphi=\varphi_0+\varphi_2+\varphi_4\in
\ka_\IC^{2*}(M)=\ka^0_\IC(M)\oplus\ka^2_\IC(M)\oplus\ka_\IC^4(M)$
satisfying:
$${\rm i)}~ d\varphi=0,~~{\rm  ii)}~ \langle\varphi,\varphi\rangle:=\varphi_2\wedge
\varphi_2-2\varphi_0\varphi_4\equiv0,{\rm ~~~and~~~} {\rm iii)}~
\langle\varphi,\bar\varphi\rangle>0.$$
\end{definition}

The inequality in iii) means that the real four-form
$\langle\varphi,\bar\varphi\rangle$ is in any point of $M$ a positive
multiple of the fixed volume form.

\begin{example}
The example that will interest us most is provided by B-field
shifts of $\sigma$. For a given ordinary
Calabi--Yau form $\sigma$ and any
real two-form $B$ the form $\exp(B)\sigma=\sigma+B\wedge\sigma$
is a generalized Calabi--Yau form. It is different from $\sigma$
only if $B^{0,2}\ne0$.
\end{example}

The quadratic form $\langle~~,~~\rangle$ defined in ii) carries
over to cohomology and yields the Mukai lattice:

\begin{definition}
Let $X$ be a K3 surface. Then the \emph{Mukai lattice}
$\widetilde H(X,\IZ)$ of $X$ is the integral cohomology
$H^*(X,\IZ)$ endowed with the quadratic form
$$\langle\varphi,\psi\rangle:=
-\varphi_0\wedge\psi_4+\varphi_2\wedge\psi_2-\varphi_4\wedge\psi_0.$$
\end{definition}

Of course, the complex structure did not matter for the definition
of the Mukai lattice, so if the K3 surface $X$ is viewed as a
complex structure on $M$, then $\widetilde
H(X,\IZ)\cong\widetilde H(M,\IZ)$. Later we shall often mean by
$\widetilde H(X,\IZ)$ the Mukai lattice of $X$ together with its
natural weight-two Hodge structure, which will be defined shortly.

Since the odd cohomology of a K3 surface is trivial, one
has
$$\widetilde H(M,\IZ)=H^2(M,\IZ)\oplus -\left(H^0(M,\IZ)\oplus
H^4(M,\IZ)\right).$$
So, as an abstract lattice it can be described by (use   $-U\cong U$):
$$\widetilde H(M,\IZ)\cong 4U\oplus 2(-E_8).$$

In analogy to the classical situation, the \emph{period of a
generalized Calabi--Yau structure} $\varphi$ on $M$ is its
cohomology class $[\varphi]$ or rather the line spanned by it
viewed as a point in $\IP(\widetilde H(M,\IC))$.

Moreover, the period of $\varphi$ can be used to introduce a Hodge
structure of weight two on $\widetilde H(M,\IZ)$, which shall be
denoted $\widetilde H(M,\varphi,\IZ)$. One defines
$$\widetilde H^{2,0}(M,\varphi):=\IC[\varphi]\subset
\widetilde H(M,\IC)$$
and then $\widetilde H^{0,2}(M,\varphi)$ is necessarily
spanned by $[\bar\varphi]$. By definition
$\widetilde H^{1,1}(M,\varphi)$ is
the orthogonal (with respect to the Mukai pairing)
complement of $\widetilde H^{2,0}(M,\varphi) \oplus\widetilde
H^{0,2}(M,\varphi)$.

\begin{examples}
i) In the case of a classical Calabi--Yau form $\varphi=\sigma$
defining a K3 surface $X$ one recovers Mukai's original
definition of the weight two Hodge structure $\widetilde H(X,\IZ)$
on the Mukai lattice, whose $(2,0)$-part is spanned by $\sigma$
and whose $(1,1)$-part is $H^{1,1}(X)\oplus (H^0\oplus H^4)(X)$.

ii) For the B-field twist of an ordinary Calabi--Yau structure
we introduce  the Hodge structure $$\widetilde H(X,B,\IZ):=
\widetilde H(X,\varphi:=\exp(B)\sigma,\IZ).$$
\end{examples}

\begin{remark}
Similar to the classical Global Torelli Theorem and eventually by
reducing to it, one proves a `generalized' Global Torelli Theorem
(see \cite{HuyGen}). This can be phrased in terms of the period
map as follows. Consider
$$\xymatrix{\widetilde\kp:
\{\IC\varphi~|~\varphi={\rm generalized~CY~form}\}/\langle{\rm
Diff}_*(M),\exp(B)\rangle\ar[r]&\widetilde Q,}$$ where $\widetilde
Q\subset \IP(\widetilde H(M,\IC))$ is the period domain
$\{x~|~\langle x,x\rangle=0,~\langle x,\bar x\rangle>0\}$,
$\widetilde\kp(\IC\varphi)=\IC[\varphi]$, and $B$ runs through
all real exact two-forms. Then $\widetilde\kp$ restricted to the
subset of those $\varphi$ satisfying a generalized K\"ahler
condition (see \cite[Sect.\ 3]{HuyGen}) is generically injective.
In analogy to Remark \ref{SurjPK3} one proves also that
$\widetilde \kp$ is surjective.

The most fascinating aspect of Hitchin's notion of generalized
Calabi--Yau structures is the occurence of classical Calabi--Yau
forms $\sigma$ as well as of symplectic generalized Calabi--Yau
forms $\exp(i\omega)$ (with $\omega$ a symplectic form) in the
same moduli space. This allows one to pass from the symplectic to
the complex world in a continuous fashion.

Note that the analogue of Siu's result has not been proved. For the time
being we do not know whether any generalized Calabi--Yau structure is
K\"ahler.
\end{remark}

\begin{remark}\label{orientRem}
The Mukai lattice $\widetilde\Gamma:=\widetilde H(M,\IZ)$ has four
positive directions and twenty negative ones. Suppose $\widetilde
\Gamma\otimes\IR\cong V_1\oplus W_1\cong V_2\oplus W_2$ are two
orthogonal decompositions of the real vector space
$\widetilde\Gamma\otimes\IR$ such that $V_1,V_2$ are
positive-definite and $W_1,W_2$ are negative-definite. Then
orthogonal projection yields an isomorphism $V_1\cong V_2$. This
allows us to compare orientations of the four-dimensional real
vector spaces $V_1$ and $V_2$. By definition, an
\emph{orientation of the positive directions} (or simply an
orientation) of the Mukai lattice is given by an orientation of a
positive-definite four-space $V\subset
\widetilde\Gamma\otimes\IR$ and two such orientations given by
orientations of $V_1,V_2\subset\widetilde\Gamma\otimes\IR$ are
equal if they correspond to each other under $V_1\cong V_2$.

If $X$ is a K3 surface, then $\widetilde H(X,\IZ)$ is naturally
endowed with an orientation (of the four positive directions).
Indeed, if $[\omega]$ is a K\"ahler class, then $\langle{\rm
Re}(\sigma),{\rm Im}(\sigma),1-\omega^2/2,\omega\rangle$ is a
positive four-dimensional subspace of $\widetilde H(X,\IR)$ and
the chosen (orthogonal) basis determines an orientation. Any other
choice of $\sigma$ or of the K\"ahler class $[\omega]$  yields
the same orientation in the above sense. Similarly, $\widetilde
H(X,B,\IZ)$ can be endowed with the orientation obtained as the
image of the previous one under $\exp(B)$.
\end{remark}
\section{Twisted K3 surfaces}\label{SecTwistedK3}

Although most of the things we shall explain or recall hold true
for a fairly general class of complex manifolds, we will restrict
to K3 surfaces. So, as before $X$ will denote a K3 surface.

 An \emph{Azumaya algebra} on $X$ is an associative,
$\ko_X$-algebra $\ka$ such that locally (in the analytic topology)
it is isomorphic to a matrix algebra ${\rm M}_r(\ko_X)$ for some
$r$. In particular, ${\ka}$ is locally free of constant rank
$r^2$.

Two Azumaya algebras are isomorphic if they are isomorphic as
$\ko_X$-algebras. By the Skolem-Noether theorem ${\rm Aut}({\rm
M}_r(\IC))\cong{\rm PGL}_r(\IC)$. Hence, isomorphism classes of
Azumaya algebras of rank $r^2$ are parametrized by the set
$H^1(X,{\rm PGL}_r)$.

To any vector bundle $E$ of rank $r$ one associates the `trivial'
Azumaya algebra ${\ka}={\cal E}nd(E)$ of rank $r^2$. This gives
rise to the following notion of equivalence between Azumaya
algebras: Two Azumaya algebras $\ka_1$ and $\ka_2$ are called
\emph{equivalent} if there exist vector bundles $E_1$ and $E_2$
such that ${\ka_1}\otimes{\ke}nd(E_1)$ and
$\ka_2\otimes{\ke}nd(E_2)$ are isomorphic Azumaya algebras.

\begin{definition}
The \emph{Brauer group} $\Br(X)$ is the set of isomorphism class
of Azumaya algebras modulo the above equivalence relation.
\end{definition}

The group structure of $\Br(X)$ is given by the tensor product of
Azumaya algebras.

A cohomological approach to Brauer groups is provided by the long
exact cohomology sequence of
$$\xymatrix{1\ar[r]&\ko_X^*\ar[r]&{\rm GL}_r\ar[r]&
{\rm PGL}_r\ar[r]& 1,}$$ which yields natural maps
$\delta_r:H^1(X,{\rm PGL}_r)\to H^2(X,\ko_X^*)$ and eventually an
injection
$$\xymatrix@W=10pt{\delta:\Br(X)\ar@{^{(}->}[r]&H^2(X,\ko_X^*).}$$
Using the commutative diagram
$$\xymatrix@H=12pt{1\ar[r]&\mu_r\ar@{^{(}->}[d]\ar[r]&{\rm
SL}_r\ar@{^{(}->}[d]\ar[r]&{\rm PGL}_r\ar[d]^=\ar[r]&1\\
1\ar[r]&\ko_X^*\ar[r]&{\rm GL}_r\ar[r]&{\rm PGL}_r\ar[r]&1}$$ one
finds that the image of $\delta_r$ is contained in the
$r$-torsion part of $H^2(X,\ko_X^*)$. Hence $\Br(X)\subset
H^2(X,\ko_X^*)$ is contained in the subgroup $H^2(X,\ko_X^*)_{\rm
tor}$  of torsion classes.

\begin{theorem}
Let $X$ be a K3 surface. Then $\Br(X)=H^2(X,\ko_X^*)_{\rm tor}$.
\end{theorem}

This result goes back to Grothendieck for projective K3 surfaces
(see \cite{GB}) and was proved in \cite{HuySch} for arbitrary K3
surfaces.

\begin{remark}
If $X$ is smooth projective or, more generally, any regular
scheme, one defines analogously the algebraic Brauer group and
compares it with the \'etale cohomology $H^2(X,\IG_m)$, which is
sometimes called the cohomological Brauer group $\Br'(X)$. The
latter contains only torsion classes and Grothendieck asked
whether the natural injection $\Br(X)\to H^2(X,\IG_m)$ is
bijective. (Without any regularity one defines $\Br'(X)$ as the
torsion part of $H^2(X,\IG_m)$.) An affirmative answer to this
question has recently been published by de Jong in \cite{dJ} for
the case of quasi-projective schemes and had earlier been proved
by Gabber (unpublished). We also recommend \cite{Lieblich2}
\end{remark}

\begin{definition}\label{DeftwistedK3}
A \emph{twisted K3 surface} $(X,\alpha)$ consists of a K3 surface
$X$ together with a class $\alpha\in H^2(X,\ko_X^*)$. We say that
$(X,\alpha)\cong (X',\alpha')$ if there exists an isomorphism
$f:X\cong X'$ with $f^*\alpha'=\alpha$.
\end{definition}

Geometrically passing from ordinary K3 surfaces to twisted K3
surfaces means that we pass from ordinary (coherent,
quasi-coherent) sheaves to twisted (coherent, quasi-coherent)
sheaves. This notion in its various incarnations
will be explained next.

\medskip

{\bf 1. Twisted sheaves.} Suppose we represent a class $\alpha\in
H^2(X,\ko_X^*)$ by a \v{C}ech 2-cocycle
$\{\alpha_{ijk}\in\Gamma(U_{ijk},\mathcal{O}^*_X)\}$ with respect
to  an open analytic covering $X=\bigcup U_i$. (Here and in the
sequel we write $U_{ij}$ and $U_{ijk}$ for the intersections
$U_i\cap U_j$ and $U_i\cap U_j\cap U_k$, respectively.)

 An \emph{$\{\alpha_{ijk}\}$-twisted
(coherent) sheaf} $E$ consists of pairs
$(\{E_i\},\{\varphi_{ij}\})$ such that the ${E}_i$ are (coherent)
sheaves on $U_i$ and
$\varphi_{ij}:{E}_j|_{U_{ij}}\rightarrow{E}_i|_{U_{ij}}$ are
isomorphisms satisfying the following conditions:
$$ {\rm i)}~\varphi_{ii}=\mathrm{id},~~ {\rm ii)}~
\varphi_{ji}=\varphi_{ij}^{-1},~~~{\rm and}~~~ {\rm iii)}~
\varphi_{ij}\circ\varphi_{jk}\circ\varphi_{ki}=\alpha_{ijk}\cdot
\mathrm{id}.$$

Morphisms between $\{\alpha_{ijk}\}$-twisted sheaves are defined
in the obvious way and one verifies that kernel and cokernel
exist. Thus, we can speak about the abelian category of
$\{\alpha_{ijk}\}$-twisted sheaves which we shall denote
$$\Coh(X,\{\alpha_{ijk}\}).$$

If $\{\alpha'_{ijk}\}$ is another \v{C}ech 2-cocycle based on the
same open covering and representing the same class $\alpha$, then
there exist $\{\lambda_{ij}\in\ko_X^*(U_{ij})\}$ such that
$\alpha'_{ijk}\alpha_{ijk}^{-1}=\lambda_{ij}\cdot\lambda_{jk}\cdot\lambda_{ki}$.
The $\{\alpha_{ijk}\}$-twisted sheaves  are in bijection with the
$\{\alpha'_{ijk}\}$-twisted sheaves via
$$\xymatrix{(\{E_i\},\{\varphi_{ij}\})\ar@{|->}[r]&(\{E_i\},\{\varphi_{ij}\cdot\lambda_{ij}\}).}$$
In particular, this yields an equivalence of abelian categories
\begin{equation}\label {secequ}
\Coh(X,\{\alpha_{ijk}\})\cong \Coh(X,\{\alpha'_{ijk}\}),
\end{equation}
which is non-canonical as it depends on the choice of
$\{\lambda_{ij}\}$.

E.g.\ if $\{\lambda_{ij}\}$ satisfies
$\lambda_{ij}\cdot\lambda_{jk}\cdot\lambda_{ki}=1$ and thus
defining a line bundle $\kl$ on $X$, then the induced equivalence
$\Coh(X,\{\alpha_{ijk}\})\cong\Coh(X,\{\alpha_{ijk}\})$ in
(\ref{secequ}) is given by the tensor product
$$\xymatrix{E\ar@{|->}[r]&
 E\otimes \kl.}$$

Similarly, if one passes to a finer open covering, then twisted
sheaves are in a natural bijection. This allows us to speak of
$\alpha$-twisted sheaves without mentioning an explicit \v{C}ech
representative of the cohomology class $\alpha\in H^2(X,\ko_X^*)$.
More precisely, all the abelian categories of twisted sheaves
with respect to some \v{C}ech cocycle representing a fixed class
$\alpha$ are equivalent, though not naturally. By abuse of
notation, we
 call the equivalence class of  these categories $\Coh(X,\alpha)$.

\medskip

{\bf 2. $\ka$-modules.} Fix $\alpha\in \Br(X)$ and pick
 a locally free coherent $\alpha$-twisted sheaf $G$
(which always exists, see below). Then
$\ka_G:=G\otimes G^*$ is an Azumaya algebra whose Brauer class is
$\alpha$. If $E$ is any $\alpha$-twisted sheaf (with respect to the
same choice of the cycle representing $\alpha$), then $E\otimes G^*$ is
an untwisted sheaf. Moreover, $E\otimes G^*$ has the structure
of an $\ka_G$-module.

The map $E\mapsto E\otimes G^*$ then defines a bijective
correspondence between $\alpha$-twisted sheaves and $\ka_G$-modules.

If we let $\Coh(X,\ka_G) $ be the abelian category of coherent
$\ka$-modules, then this map yields an equivalence
$$\Coh(X,\ka_G)\cong\Coh(X,\alpha).$$
\medskip

{\bf 3. Sheaves on gerbes.} To an Azumaya algebra $\ka$ as well as
to a cocycle $\{\alpha_{ijk}\}$ representing a class $\alpha\in
H^2(X,\ko_X^*)$ one can associate $\IG_m$-gerbes over $X$, which
are called $\km_\ka$ and $\km_{\{\alpha_{ijk}\}}$, respectively.

The gerbe $\km_\ka$ associates to $T\to X$  is the category
$\km_\ka(T)$ whose objects are pairs $(E,\eta)$ with $E$ a locally
free coherent sheaf on $T$ and $\eta:\ke nd(E)\cong\ka_T$ an
isomorphism of $\ko_T$-algebras (see \cite{Giraud,Milne}). A
morphism $(E,\eta)\to (E',\eta')$ is given by an isomorphism
$E\cong E'$ that commutes with the $\ka_T$-actions induced by
$\eta$ and $\eta'$. It is easy to see that the group of
automorphisms of an object $(E,\eta)$ is $\ko^*(T)$.

The gerbe $\km_{\{\alpha_{ijk}\}}$ associates to $T\to X$ the
category $\km_{\{\alpha_{ijk}\}}(T)$ whose objects are collections
$\{\kl_i,\varphi_{ij}\}$ with $\kl_i\in\Pic(T_{U_i})$ and
$\varphi_{ij}:\kl_j|_{T_{U_{ij}}}\cong\kl_i|_{T_{U_{ij}}}$
satisfying
$\varphi_{ij}\cdot\varphi_{jk}\cdot\varphi_{ki}=\alpha_{ijk}$
(see \cite{Lieblich}). A morphism
$\{\kl_i,\varphi_{ij}\}\to\{\kl'_i,\varphi'_{ij}\}$ is given by
isomorphisms $\kl_i\cong\kl_i'$ compatible with $\varphi_{ij}$
and $\varphi_{ij}'$. For yet another construction of a gerbe
associated to $\alpha$ see \cite{dJ}.

 Any sheaf $\kf$ on a
$\IG_m$-gerbe $\km\to X$ comes with a natural $\IG_m$-action and
thus decomposes as $\kf=\bigoplus\kf^m$, where the $\IG_m$-action
on $\kf^m$ is given by the character $\lambda\mapsto\lambda^m$. If
$\kf=\kf^m$, then $\kf$ is called of weight $m$.

There are  natural bijections (of sets of isomorphism classes)
$$\xymatrix{\{\ka-{\rm modules}\}\ar@{<->}[r]&\{{\rm sheaves~on~}
\km_\ka~{\rm of~weight~one}\}}$$
$$\xymatrix{\{\{\alpha_{ijk}\}-{\rm twisted~sheaves}\}\ar@{<->}[r]&
\{{\rm sheaves~on~}\km_{\{\alpha_{ijk}\}}~{\rm of~weight~one}\},}$$ which
hold  for (quasi)-coherent sheaves. This yields
equivalences $$\Coh(X,\ka)\cong\Coh(\km_\ka)_1~~~~~~{\rm
and}~~~~~~
\Coh(X,\{\alpha_{ijk}\})\cong\Coh(\km_{\{\alpha_{ijk}\}})_1.$$
Here, the abelian categories on the gerbes are the categories of
coherent sheaves of weight one. One can also show that
$\Coh(X,\{\alpha^\ell_{ijk}\})\cong\Coh(\km_{\{\alpha_{ijk}\}})_\ell$.
\cite{dJ,DP,Lieblich} for more details.

 Isomorphism classes of
$\IG_m$-gerbes are in bijection with classes in $H^2(X,\ko^*_X)$
(see \cite{Milne}) and the isomorphism class of $\km_\ka$
corresponds indeed to $\delta[\ka]$. Similarly, the above
construction of $\km_{\{\alpha_{ijk}\}}$ ensures that its
isomorphism class corresponds to the class $\alpha$.

In order to construct
a concrete gerbe in the isomorphism class determined by $\alpha$
we had to choose a specific Azumaya algebra $\ka$ or
 a cocycle $\{\alpha_{ijk}\}$ representing $\alpha$.
In the first case we have to assume $\alpha\in\Br(X)$.

One can construct directly an isomorphism
$\km_\ka\cong\km_{\{\alpha_{ijk}\}}$ for
 an appropriate  cocycle $\alpha_{ijk}$ representing
$\alpha=\delta[\ka]$. This goes as follows: Choose a covering
$X=\bigcup U_i$, isomorphisms $\eta_i:\ke nd(E_i)\cong\ka_{U_i}$,
where the $E_i$ are (locally) free, and isomorphisms
$\xi_{ij}:E_j|_{U_{ij}}\cong E_i|_{U_{ij}}$ compatible with the
$\eta_i$. Then $\xi_{ij}\cdot\xi_{jk}\cdot\xi_{ki}$ is given by
multiplication with a scalar function $\alpha_{ijk}$ whose
cohomology class represents $\delta[\ka]$. If
$(E,\eta)\in\km_\ka(T)$, then $E_i\cong E_{U_i}\otimes\kl_i$ for
certain line bundles $\kl_i$ on $T_{U_i}$ and one may choose
isomorphisms $\varphi_{ij}$ between them (over $U_{ij}$) such
that $\xi_{ij}$ is ${\rm id}\otimes \varphi_{ij}$. Associating to
$(E,\eta)$ the collection $\{\kl_i,\varphi_{ij}\}$ defines an
isomorphism of gerbes $\km_\ka\to\km_{\{\alpha_{ijk}\}}$.

\medskip

{\bf 4. Sheaves on the Brauer--Severi variety.} The following has
been explained in Yoshioka's article \cite{Yoshioka}. Suppose
$E=(\{E_i\},\{\varphi_{ij}\})$ is a locally free
$\{\alpha_{ijk}\}$-twisted sheaf. The projective bundles
$\IP(E_i)\to U_i$ glue to the Brauer--Severi variety
$\pi:\IP(E)\to X$ and the relative tautological line bundles
$\ko_{\IP(E_i)}(1)$ glue to a
$\{\pi^*\alpha^{-1}_{ijk}\}$-twisted line bundle $\ko_\pi(1)$ on
$\IP(E)$. If $F=(\{F_i\},\{\psi_{ij}\})$ is any
$\{\alpha_{ijk}\}$-twisted sheaf, then $\pi^*F\otimes\ko_\pi(1)$
is a true sheaf in a natural way. This yields an equivalence of
$\Coh(X,\{\alpha_{ijk}\})$ with the full subcategory
$\Coh(\IP(E)/X)$ of $\Coh(\IP(E))$ of all coherent sheaves $F'$
on $\IP(E)$ for which the natural morphism
$\pi^*\pi_*(F'\otimes(\pi^*E\otimes\ko_\pi(1))^*)\to
F'\otimes(\pi^*E\otimes\ko_\pi(1))^*$ is an isomorphism:
$$\Coh(X,\{\alpha_{ijk}\})\cong\Coh(\IP(E)/X).$$
Note that the bundle $\pi^*E\otimes\ko_\pi(1)$ can be described
as the unique non-trivial extension $0\to\ko_{\IP(E)}\to
\pi^*E\otimes\ko_\pi(1)\to\kt_{\IP(E)}\to 0$ and thus depends
only on the Brauer--Severi variety $\pi:\IP(E)\to X$.

Isomorphism classes of $\IP^r$-bundles are also parametrized by
$H^1(X,{\rm PGL}_r)$. Locally a $\IP^r$-bundle is of the form
$\IP(E_i)$ for some locally free sheaf $E_i$. They glue to a
$\delta(\IP)$-twisted sheaf.

\section{Twisted Chern characters}\label{SectTwistChern}

In order to study twisted K3 surfaces and twisted sheaves by
cohomological methods, one needs a good cohomology theory and the
notion  of twisted Chern characters.

Let us begin by introducing the weight two Hodge structure
$\widetilde H(X,\alpha,\IZ)$ of a twisted K3 surface
$(X,\alpha)$. The exponential sequence shows that any element
$\alpha\in H^2(X,\ko_X^*)$ can be written as $\exp(B^{0,2})$ for
some B-field $B\in H^2(X,\IR)$. If $\alpha$ is a torsion class,
we may choose $B$ to be rational. The exponential sequence also
shows that a given $B$ may be changed by integral B-fields
$B_0\in H^2(X,\IZ)$ without changing the Brauer class $\alpha$.
Once a B-field lift $B\in H^2(X,\IR)$ of a class $\alpha$ is
chosen, one considers the generalized Calabi--Yau form
$\exp(B)\sigma$ and its natural weight two Hodge structure
$\widetilde H(X,B,\IZ)$. For $B_0\in H^2(X,\IZ)$ multiplication
with the integral class $\exp(B_0)$ defines a Hodge isometry
$$\widetilde H(X,B,\IZ)\cong \widetilde H(X,B+B_0,\IZ).$$ This
allows us to introduce
$$\widetilde H(X,\alpha,\IZ)$$
as the Hodge isometry type of $\widetilde H(X,B,\IZ)$ with $B$ an
arbitrary B-field lift of $\alpha$. We emphasize that this is an
abstract Hodge structure, for the realization of which one needs
to choose a concrete B-field lift of $\alpha$.

There are various approaches towards twisted Chern characters,
e.g.\ \cite{Heinl,HS1,Yoshioka}. The one introduced in \cite{HS1}
seems not very canonical, as it depends on the additional choice
of a B-field. It is, however, the one that works best in the
context of twisted K3 surfaces, as it allows us to work with
integral(!) Hodge structures.

Let $B\in H^2(X,\IQ)$ and $\alpha\in H^2(X,\ko_X^*)$ be the
induced Brauer class, i.e.\ the image of $B^{0,2}\in
H^2(X,\ko_X)$ under the exponential $H^2(X,\ko_X)\to
H^2(X,\ko_X^*)$. Equivalently, $\alpha$ is the image of $B\in
H^2(X,\IQ)$ under the composition of the exponential map
$H^2(X,\IQ)\to H^2(X,\IC^*)$ and the natural inclusion
$\IC^*\subset\ko_X^*$. In addition, choose a \v{C}ech cocycle
$B_{ijk}\in \Gamma(U_{ijk},\IQ)$ representing $B$ and let
$\alpha_{ijk}:=\exp(B_{ijk})$ be the induced \v{C}ech cocycle
representing $\alpha$.

Once this $\{\alpha_{ijk}\}$ is fixed, we can speak of
$\{\alpha_{ijk}\}$-twisted sheaves and we aim at defining a
twisted Chern character for those. Before we can do this in
practice we need to make yet another choice. Viewing the
$B_{ijk}$ as differentiable functions allows us to write them as
$B_{ijk}=-a_{ij}+a_{ik}-a_{jk}$ for certain differentiable
functions $a_{ij}:U_{ij}\to\IR$. (We use $H^2(X,\kc^\infty_X)=0$
and might have to refine the covering.)

Let now $E=(\{E_i\},\{\varphi_{ij}\})$ be an
$\{\alpha_{ijk}\}$-twisted sheaf. Then one defines
$$E_B:=(\{E_i\},\{\varphi'_{ij}:=\varphi_{ij}\cdot\exp(a_{ij})\})
~~~{\rm and}~~{\rm ch}^B(E):={\rm ch}(E_B).$$ Note that $E_B$
describes an untwisted sheaf, for the $\varphi'_{ij}$ satisfy the
usual cocycle condition. First observe that this definition of
the twisted Chern character is independent of the choice of
$\{a_{ij}\}$. Indeed, passing to $a_{ij}+a_{ij}'$ with
$-a'_{ij}+a'_{ik}-a'_{jk}=0$ would change the bundle $E_B$ by a
twist with the line bundle $L$ corresponding to the cocycle
$\{\exp(a_{ij}')\}$. Since ${\rm
c}_1(L)=\{-a'_{ij}+a'_{ik}-a'_{jk}\}=0$, this has no effect on
${\rm ch}(E_B)$.

On the other hand, although not reflected by our notation, the
choice of $\{B_{ijk}\}$ and the resulting $\{\alpha_{ijk}\}$ is
necessary in order to be able to define ${\rm ch}^B(E)$, as we
could otherwise not speak about $\{\alpha_{ijk}\}$-twisted sheaf.
However, as was explained earlier, there is a non-canonical
bijection between $\{\alpha_{ijk}\}$-twisted sheaves and
$\{\alpha'_{ijk}\}$-twisted sheaves for two \v{C}ech cocycles
representing $\alpha$ and  and our Chern character ${\rm ch}^B$
is compatible with it.

Indeed, if $B_{ijk}':=B_{ijk}+(b_{ij}-b_{ik}+b_{jk})$, then
$\{\alpha_{ijk}\}$ and $\{\alpha'_{ijk}\}$ differ by the boundary
of $\{\lambda_{ij}:=\exp(b_{ij})\}$ and we may send
$E=(\{E_i\},\{\varphi_{ij}\})$ to the $\{\alpha'_{ijk}\}$-twisted
sheaf $E'=(\{E_i\},\{\varphi_{ij}\cdot\lambda_{ij}\})$. (The
given modification of $\{B_{ijk}\}$ by the boundary of
$\{b_{ij}\}$ induces a canonical bijection between
$\{\alpha_{ijk}\}$-twisted sheaves and
$\{\alpha'_{ijk}\}$-twisted sheaves, which otherwise does not
exist.) Clearly, $E_B$ and $E'_B$ are defined by the same
cocycle. Thus, ${\rm ch}^B$ does not depend on the \v{C}ech
cocycle representing $B$.

The following properties of the twisted Chern character ${\rm
ch}^B$ have been observed in \cite{HS1}:

{\rm i)} ${\rm ch}^B(E_1\oplus E_2)={\rm ch}^B(E_1)+{\rm
ch}^B(E_2)$.

{\rm ii)} If $B={\rm c}_1(L)\in H^2(X,\IZ)$, then ${\rm
ch}^B(E)=\exp({\rm c}_1(L))\cdot {\rm ch}(E)$.

{\rm iii)} ${\rm ch}^{B_1}(E_1)\cdot {\rm ch}^{B_2}(E_2)={\rm
ch}^{B_1+B_2}(E_1\otimes E_2)$.

\begin{remark}
i) In the note \cite{Heinl} Heinloth explains the relation between
the twisted Chern character ${\rm ch}^B$ and the usual Chern
character on the gerbe. He first proves that
\begin{equation}\label{Heinlothequ}
H^*(\km_\alpha,\IQ)\cong H^*(X,\IQ)[z],
\end{equation}
where $z={\rm c}_1(E)$ with $E$ some vector bundle of weight one
on $\km_\alpha$. In particular, the isomorphism in
(\ref{Heinlothequ}) depends on this choice. He furthemore
explains that the choice of a differentiable line bundle $L$ of
weight one on $\km_\alpha$ allows to define ${\rm ch}_L(E)$ as
${\rm ch}(E\otimes L^*)$, which makes sense as $E\otimes L^*$ has
weight zero, i.e.\ comes from $X$. The choice of $L$ corresponds
to the choice of the B-field $B$ and one obtains ${\rm
ch}^B(E)={\rm ch}_L(E)$.

ii) Yoshioka uses yet other conventions to descend from the derived
category to cohomology. A detailed comparison of the various twisted
Chern characters can be found in \cite{HS2}.
\end{remark}

We are primarily interested in twisted K3 surfaces $(X,\alpha)$
and their natural weight two Hodge structure $\widetilde
H(X,\alpha,\IZ)$. Unfortunately, I don't know of any elegant way
to define directly ${\rm ch}^\alpha:\Coh(X,\alpha)\to \widetilde
H(X,\alpha,\IZ)$. Twisted Chern character and twisted cohomology
can physically only be realized after choosing a B-field lift.

\bigskip

The twisted as well as the untwisted Chern character has to be
modified to the Mukai vector. Only working with the Mukai vector
allows one to descend from equivalences of derived categories to
isomorphisms of cohomologies.

With the above notation one defines the \emph{Mukai vector} as
$$v^B(~~):={\rm ch}^B(~~)\cdot\sqrt{{\rm td}(X)}.$$
It can be applied to (twisted) sheaves as well as to complexes of
sheaves and maps coherent (twisted) sheaves to classes in
$\widetilde H^{1,1}(X,B,\IZ)$. The
definition in the untwisted case, i.e.\ $B=0$, is due to Mukai
and we write simply $v(~~)$ in this case. Note that for K3
surfaces $\sqrt{{\rm td}(X)}=(1,0,1)$.

If $E,F\in\Coh(X,\alpha)$, then the Hirzebruch--Riemann--Roch
formula reads
$$\chi(E,F):=\sum(-1)^i\dim\Ext^i(E,F)=-\langle v^B(E),v^B(F)\rangle.$$

\begin{examples}
The following (untwisted) Mukai vectors of sheaves on K3 surfaces
are used frequently: i) $v(\ko_X)=(1,0,1)$, ii)
$v(\ko_C(-1))=(0,[C],0)$, iii) $v(k(x))=(0,0,1)$.
\end{examples}
\section{The bounded derived category of a K3 surface}

Let $X$ be a K3 surface and let $\alpha\in H^2(X,\ko_X^*)$ be a
class represented by a \v{C}ech cocycle $\{\alpha_{ijk}\}$. One
associates to $X$ and $(X,\alpha)$ the abelian categories
$\Coh(X)$ and $\Coh(X,\{\alpha_{ijk}\})$ of coherent sheaves on
$X$ and $\{\alpha_{ijk}\}$-twisted coherent sheaves, respectively.

\begin{remark}
These categories are defined for any $X$ and any $\alpha$. However, for $X$ not
projective or $\alpha$ not torsion, there exist usually very few interesting (twisted)
coherent sheaves on $X$. The categories $\Coh(X)$ and $\Coh(X,\alpha)$ will simply be
too small to be interesting. Thus, studying these categories is a sensible thing
to do only under these `projectivity' assumptions.
\end{remark}

The following theorem in the untwisted case is a special case of
a result of Gabriel \cite{Gabriel}.

\begin{theorem} Suppose $(X,\alpha)$ and $(X',\alpha')$ are two twisted K3 surfaces with
$X,X'$ projective and $\alpha,\alpha'$ torsion. Then $(X,\alpha)\cong(X',\alpha')$
if and only if there exists an equivalence $\Coh(X,\alpha)\cong\Coh(X',\alpha')$.
\end{theorem}

The basic idea is very simple. Consider the \emph{minimal
objects} in $ \Coh(X,\alpha)$, i.e.\ those objects that do not
contain any proper non-trivial sub-objects. This is certainly a
notion that is preserved under any equivalence. On the other hand,
it is straightforward to prove that minimal objects in
$\Coh(X,\alpha)$ are just the skyscraper sheaves $k(x)$ of closed
points $x\in X$. Hence, any equivalence will induce a bijection
$X\cong X'$. It remains to show that this natural bijection is a
morphism. In the case of surfaces the topology is determined by
points and curves. Since curves have trivial Brauer group, every
curves supports an invertible twisted sheaf.

\begin{remark}
i) The result holds in much broader generality. Gabriel proves
the untwisted case for arbitrary schemes. For the twisted case
see \cite{Perego}.

ii) If $X$ is not projective, then the abelian category $\Coh(X)$
does not, in general, encode the variety completely. E.g.\ as
shown by Verbitsky in \cite{Verb} two very general complex tori
 have equivalent abelian categories.
\end{remark}

The upshot is that passing from (twisted) K3 surfaces to their abelian categories
no information is lost. We still can formulate a Global Torelli Theorem:

\begin{corollary}
Suppose $(X,\alpha)$ and $(X',\alpha')$ are two twisted K3
surfaces with $X,X'$ projective and $\alpha,\alpha'$ torsion
classes. Then $\Coh(X,\alpha)\cong\Coh(X',\alpha')$ if and only if
there exists a Hodge isometry $H^2(X,\IZ)\cong H^2(X',\IZ)$ such
that the induced map $\Br(X)\to \Br(X')$ sends $\alpha$ to
$\alpha'$.
\end{corollary}

One can be a little more specific: Any equivalence
$\Coh(X,\alpha)\cong\Coh(X',\alpha')$ is of the form
$(M\otimes~~~~)\circ f_*$,
where $f:X\cong X'$ is an isomorphism and $M\in \Pic(X')$.
Note that a priori it is not clear how to pass from an
equivalence of the abelian categories
directly to a Hodge isometry. In fact, the natural
Hodge isometry associated to
an equivalence of the form  $(M\otimes~~~~)\circ f_*$
would not be $f_*$, but
$\exp({\rm c}_1(M))\circ f_*$. Here, $\exp({\rm c}_1(M))$ means multiplication with
the Chern character ${\rm ch}(M)$.

\begin{remark}
We have chosen to work with  $\{\alpha_{ijk}\}$-twisted sheaves,
but rephrasing everything in terms of $\ka$-modules or sheaves on
the gerbes $\km_\ka$ or $\km_{\{\alpha_{ijk}\}}$ or on the
Brauer--Severi variety $\IP(E)$ is possible. One would then
consider the abelian categories $\Coh(X,\ka)$, $\Coh(\km_\ka)_1$,
$\Coh(\km_{\{\alpha_{ijk}\}})_1$, or $\Coh(\IP(E)/X)$ (see Section
\ref{SecTwistedK3}). It follows from the discussion there that
all these categories are equivalent if $\delta[\ka]=\alpha$ and
$E\in\Coh(X,\{\alpha_{ijk}\})$. It is largely a matter of taste
which one is prefered.

Let us emphasize, however, that there is no $\IG_m$-gerbe
$\km_\alpha$ naturally associated to a Brauer class $\alpha\in
\Br(X)$ but only an isomorphism class of $\IG_m$-gerbes. In
particular, before being able to introduce the abelian category
of the twisted K3 surface $(X,\alpha)$ one has to make a choice,
either of a cocycle $\{\alpha_{ijk}\}$ representing $\alpha$, of
an Azumaya algebra $\ka$ with $\delta[\ka]=\alpha$, of a
$\IG_m$-gerbe in the isomorphism class determined by $\alpha$, or
of a Brauer--Severi variety $\pi:\IP(E)\to X$ realizing $\alpha$.
\end{remark}

\bigskip

Let us  stick for the rest of this section to ordinary K3
surfaces. The twisted case will be discussed in the next section.
This is historically correct and makes, I hope, the most general
case easier to digest

\begin{definition}
The \emph{derived category} $\Db(X)$ of a K3 surface $X$ is the bounded
derived category of the abelian category $\Coh(X)$, i.e.\
$$\Db(X):=\Db(\Coh(X)).$$
\end{definition}

The category $\Db(X)$ is a $\IC$-linear triangulated category and
equivalences between such categories will always assumed to be
$\IC$-linear and exact, i.e.\ shifts and distinguished triangles are respected.
Two K3 surfaces are called \emph{derived equivalent} if there exists
an equivalence $\Db(X)\cong\Db(X')$.

Why passing from the abelian category $\Coh(X)$ to its derived
category might change things, is explained by Mukai's celebrated
example. It marked the beginning of the theory of Fourier--Mukai
transforms (see \cite{MukaiAV}): Let $A$ be an abelian variety and
let $\widehat A$ be its dual abelian variety. In general, $A$ and
$\widehat A$ are non-isomorphic. Indeed, they are isomorphic if
and only if $A$ is principally polarized. Nevertheless, there
always exists an exact equivalence
$$\Db(A)\cong\Db(\widehat A).$$
Mukai not only proves the equivalence of the two derived
categories, but suggests how to produce geometrically interesting
equivalences in general. This has led to the concept of
Fourier--Mukai transforms.

\begin{definition}
Let $X$ and $X'$ be any two smooth projective varieties and
$\kp\in\Db(X\times X')$. The \emph{Fourier--Mukai} transform with
\emph{Fourier--Mukai kernel} $\kp$ is the exact functor:
$$\xymatrix{\Phi_\ke:\Db(X)\ar[r]& \Db(X')},~~\xymatrix{E^\bullet\ar@{|->}[r]&
Rp_*(q^*E^\bullet\otimes^L\kp),}$$
where $q$ and $p$ denote the two projections from $X\times X'$.
\end{definition}

In general, a Fourier--Mukai transform will not define an
equivalence, but due to a deep theorem of Orlov the converse
holds, see \cite{Or}:

\begin{theorem}\label{OrlovExist}
Let $\Phi:\Db(X)\cong\Db(X')$ be an exact equivalence. Then there
exists a unique object $\kp\in\Db(X\times X')$ (up to
isomorphism) such that $\Phi\cong\Phi_\kp$.
\end{theorem}

Note that this is somewhat equivalent to the fact that any
equivalence between the abelian categories of coherent sheaves
has a very special form (composition of an isomorphism with a
line bundle twist). However, the object $\kp$ might be difficult
to describe explicitly and in general it will be a true complex
(and not just a shifted sheaf).

In many situations the following criterion can be used to decide
whether a given Fourier--Mukai transformd defines an equivalence.
See the original articles \cite{BO2,Bridgeland} or \cite{HuyFM}
for the proof and similar results.

\begin{theorem}\label{Bridgcrit}
Suppose the Fourier--Mukai transform $\Phi_\kp:\Db(X)\to\Db(X')$ satisfies the following two
conditions:

 {\rm i)} $\dim\Hom(\Phi(k(x)),\Phi(k(x)))=1$ for any
$x\in X$ and

 {\rm ii)} $\Hom(\Phi(k(x)),\Phi(k(y))[i])=0$ for $x\ne
y$ or $i<0$ or $i>\dim(X)$.

Then $\Phi_\kp$ defines an equivalence.
\end{theorem}

\begin{examples}\label{FMExK3}
i) Let $A$ be an abelian variety and let $\widehat A$ be its dual.
The Poincar\'e bundle $\kp$ on $\widehat A\times A$ can be
considered as an object in $\Db(\widehat A\times A)$. The famous
result of Mukai alluded to before states that the induced
Fourier--Mukai transform $\Phi_\kp:\Db(\widehat A)\to\Db(A)$ is an
equivalence. Nowadays the result can be obtained as a direct
consequence of
 Theorem \ref{Bridgcrit}.

ii) Any isomorphism $X\cong X'$ induces an equivalence
$\Db(X)\cong\Db(X')$. The Fourier--Mukai kernel is the
structure sheaf of its graph.

iii) Suppose $M$ is a moduli space of stable sheaves
on a K3 surface $X$. If $M$ is complete and two-dimensional, then
$M$ is a K3 surface (see \cite{Mu} or \cite{HL}). If $M$ is fine, i.e.\
a universal sheaf $\ke$ on $M\times X$ exists, then Theorem \ref{Bridgcrit}
again applies and yields an equivalence $\Phi_\ke:\Db(M)\cong\Db(X)$.

This is in analogy to i), where $\widehat A$ can be considered as
a moduli space of line bundles on $A$ and $\kp$  as a universal
family. For K3 surfaces however the `dual' K3 surface provided by
a moduli space $M$ as above is not unique.

iv) If $L$ is a line bundle on a projective variety, then
$F^\bullet \mapsto L\otimes F^\bullet$ defines an equivalence
$\Db(X)\cong\Db(X)$ which can be described as the Fourier--Mukai
transform with kernel $\iota_*L$, where $\iota:X\to X\times X$ is
the diagonal embedding.

v) Suppose $X$ is a K3 surface containing an irreducible smooth
rational curve $C\subset X$. Consider the tautological line bundle
$\ko_C(-1)$ on $C\cong\IP^1$ as an object in $\Db(X)$. The trace
induces a natural morphism
$\ko_C(-1)\boxtimes\ko_C(-1)\dual\to\ko_\Delta$ in $\Db(X\times
X)$, where $\ko_C(-1)\dual$ denotes the derived dual. The cone of
this morphism shall be denoted $\kp_{\ko_C(-1)}\in \Db(X\times
X)$ and the induced Fourier--Mukai functor is the \emph{spherical
twist} $T_{\ko_C(-1)}:=\Phi_{\kp_{\ko_C(-1)}}$, which is an
equivalence.

vi) The sheaf $E:=\ko_C(-1)$ in v) is a \emph{spherical object},
i.e.\ $E$ satisfies $$\Ext^*_X(E,E)=H^*(S^2,\IC).$$ Kontsevich
proposed to consider a spherical twist $T_E$ associated to any
spherical object $E$ on a Calabi--Yau manifold and Seidel and
Thomas were able to prove that $T_E$ is indeed an equivalence.
(This time Theorem \ref{Bridgcrit} is of no use, another kind of
spanning class is needed here, see \cite{ST} or \cite[Ch.\
8]{HuyFM}.) Other examples of spherical objects on a K3 surface
$X$ are $\ko_X$, or more generally any line bundle, and simple
rigid sheaves.
\end{examples}

Any Fourier--Mukai transform $\Phi_\kp:\Db(X)\to \Db(X')$ induces
the coho\-mo\-logical Fourier--Mukai transform
$$\xymatrix{\Phi_\kp^H:H^*(X,\IQ)\ar[r]& H^*(X',\IQ),}$$
which is defined in terms of the Mukai vector $v(\kp)\in
H^*(X\times X',\IQ)$ as
$$\xymatrix{\gamma\ar@{|->}[r]&p_*(v(\kp).q^*(\gamma)).}$$

If $\Phi_\kp$ is an equivalence, then $\Phi_\kp^H$ is bijective.
This innocent looking statement is not trivial, as the part in
$H^*(X,\IQ)$ that comes from objects in $\Db(X)$ may be small.
Note that in general $\Phi_\kp$ does neither respect the grading
nor the algebra structure nor will it be defined over $\IZ$.

\bigskip

Back to the case of K3 surfaces, one finds that in  Examples
\ref{FMExK3}, v) and vi) the cohomological spherical shift $T^H_E$
is given by the reflection $\gamma\mapsto
\gamma+(\gamma.v(E))v(E)$. In particular,
$T^H_{\ko_C(-1)}=s_{[C]}$. The tensor product $L\otimes~~$ acts
by multiplication with $\exp({\rm c}_1(L))$ on $H^*(X,\IZ)$.

\begin{remark}
In ii), Remark \ref{RemGT} we have pointed out that the Hodge
isometry $s_{[C]}$ is, for trivial reasons, not induced by any
automorphism of $X$. This is cured by the above observation which
says that it can be lifted, however, to an autoequivalence of
$\Db(X)$.
\end{remark}

Mukai shows in \cite{Mu} that $\Phi_\kp^H$ of a derived
equivalence $\Phi_\kp$ between two K3 surfaces is defined over
$\IZ$ and that it defines a Hodge isometry $\widetilde
H(X,\IZ)\cong \widetilde H (X',\IZ)$. Combined with Orlov's
result this becomes

\begin{corollary}\label{MukaiDerHodge}
Any derived equivalence $\Phi:\Db(X)\cong\Db(X')$ between two K3 surfaces
induces naturally a Hodge isometry
$$\xymatrix{\Phi:\widetilde H(X,\IZ)\ar[r]^-\sim& \widetilde H(X',\IZ).}$$
\end{corollary}

The other direction, namely how to deduce from
the existence of a Hodge iso\-metry
of the Mukai lattices of two K3 surfaces
the existence of a derived equivalence,
was proved by Orlov. Both results together combine to

\begin{theorem}\label{DGTT}{\bf (Derived Global Torelli)}
Two projective K3 surfaces $X$ and $X'$ are derived equivalent if
and only if there exists a Hodge isometry $$\widetilde
H(X,\IZ)\cong \widetilde H(X',\IZ).$$
\end{theorem}

For the complete proof of the theorem the reader may consult the
original article \cite{Or} or \cite[Ch.\ 10]{HuyFM}. What is
important to know for our purpose is that a given Hodge isometry
is modified by the cohomological Fourier--Mukai transforms of the
type iii)-v) in Example \ref{FMExK3} such that the new Hodge
isometry induces a Hodge isometry of the standard weight-two
Hodge structure $H^2(X,\IZ)$ of $X$ with the one of some moduli
space $Y$ of sheaves on $X'$, which is again a K3 surface. Then
the classical Global Torelli (see Theorem \ref{CGT}) applies and
yields an isomorphism $X\cong Y$.

\begin{remark}\label{Remarkderorient}
There is a minor, but annoying issue in the argument. At the very
end one has to ensure that the image of a K\"ahler class is a
K\"ahler class and not only up to sign. But of course, composing
the given Hodge isometry with the Hodge isometry ${\rm
id}_{H^0}\oplus-{\rm id}_{H^2}\oplus {\rm id}_{H^4}$ clears this
problem.

The problem is reflected by the following more precise result
that is the derived analogue of iii), Remark \ref{RemGT}: Any
orientation preserving Hodge isometry $\widetilde
H(X,\IZ)\cong\widetilde H(X',\IZ)$ lifts to a derived
equivalence, i.e.\ is of the form $\Phi^H_\kp$. (Note that the
uniqueness of the classical Global Torelli Theorem does not hold.
See the discussion in Section \ref{SecEnd}.)

In \cite{Sz} Szendr\H{o}i suggested that the mirror symmetry
analogue of the result of Donaldson (\ref{BorceaDonald}) should
say that a Hodge isometry not preserving the natural orientation
cannot be lifted. In other words, one expects that $\Phi^H$ in
Corollary \ref{MukaiDerHodge} is always orientation preserving.
That this is the case at least for all known examples was proven
in \cite{HS1}. The Fourier--Mukai equivalence induced by the
universal family of stable sheaves is the only non-trivial case.

The issue becomes more serious in the twisted case. Composing a
given Hodge isometry $g$ with $g_0:={\rm id}_{H^0}\oplus-{\rm
id}_{H^2}\oplus {\rm id}_{H^4}$ in order to reverse the
orientation is not allowed anymore. Indeed, only in the untwisted
case is $g_0$  naturally a Hodge isometry.
\end{remark}

If one prefers to work with the transcendental part
of the Hodge structure, then the above theorem becomes

\begin{corollary}\label{GTtrans}
Two projective K3 surfaces $X$ and $X'$ are derived
equivalent if and only if there exists a Hodge isometry
$T(X)\cong T(X')$.
\end{corollary}

\begin{proof}
Recall that the transcendental lattice $T$
of a Hodge structure of weight two on a lattice $H$ is
the orthogonal complement of $H\cap H^{1,1}$ and that $T$
is again a Hodge structure of weight two.
In our geometric situation, $T(X)$ is the transcendental
lattice of $H^2(X,\IZ)$ or, equivalently, of $\widetilde H(X,\IZ)$.

Clearly, any Hodge isometry $\widetilde H(X,\IZ)\cong
\widetilde H(X',\IZ)$ induces a Hodge isometry $T(X)\cong T(X')$.
Conversely, due to a result of Nikulin
 \cite{Nikulin}, any Hodge isometry $T(X)\cong T(X')$ can be extended
to a Hodge isometry of the Mukai lattices $\widetilde H(X,\IZ)\cong
\widetilde H(X',\IZ)$. The reason behind this is the existence
of the hyperbolic plane $H^0\oplus H^4$ in the orthogonal complement
of $T(X)\subset  \widetilde H(X,\IZ)$. Note that in general the
orthogonal complement of $T(X)\subset H^2(X,\IZ)$ does not contain
any hyperbolic plane, which explains why derived equivalent K3 surfaces are not necessarily isomorphic.
\end{proof}

As a consequence of the proof of the theorem, Orlov obtains

\begin{corollary}
Two projective K3 surfaces $X$ and $X'$ are derived equivalent if and only if
$X'$ is isomorphic to a moduli space of stable sheaves on $X$.
\end{corollary}

The polarization needs to be fixed appropriately and the sheaves
might a priori be torsion, e.g.\ $X$ itself is viewed as the
moduli space of skyscraper sheaves $k(x)$ or, equivalently, of
the ideal sheaves ${\cal I}_x$ of closed points $x\in X$. In
fact, in the corollary one could replace `stable sheaves' by
`torsion free stable sheaves'. \footnote{More recently, we could
show in \cite{Huyab} that this can be further improved to
`$\mu$-stable locally free'.}

\begin{remark}\label{numberFM}
Let us also mention that a general conjecture stating that any
smooth projective variety admits only finitely many Fourier--Mukai
partners, i.e.\ smooth projective varieties with equivalent
derived categories, can be proved for K3 surfaces. Once Corollary
\ref{MukaiDerHodge} is established, one uses lattice theory. A
related natural question asks for the number of isomorphism types
of K3 surfaces derived equivalent to a given K3 surface $X$ in
terms of the period of $X$. This question has been addressed in
\cite{Oguisonew,St}.
\end{remark}
\bigskip

One of the standard examples of Fourier--Mukai kernels defining a
derived equi\-valence between K3 surfaces is provided by the
universal sheaf $\ke$ on $X\times M$, where $M$ is a complete,
fine moduli space of stable sheaves of dimension two (see iii),
Examples \ref{FMExK3}).

There are examples of  moduli spaces $M$ of stable sheaves on a
K3 surface $X$ which are not fine, i.e.\ a universal sheaf $\ke$
does not exist. Locally in the analytic (or \'etale) topology of
$M$ one finds universal sheaves, but the obstruction to glue those
to a global universal sheaf might be non-trivial.
C\u{a}ld\u{a}raru observed in \cite{C} that this obstruction can
be considered as a Brauer class $\alpha\in H^2(M,\ko_M^*)$ and
that the local universal sheaves glue to a
$(1\times\alpha)$-twisted universal sheaf $\ke$ on $X\times M$.

If $M$ is complete and two-dimensional, then $M$ is a K3 surface,
but in general not derived equivalent to $X$. However, in
\cite{C} it is shown that the twisted universal sheaf $\ke$
induces a Fourier--Mukai transform that does define an
equivalence $$\Db(M,\alpha^{-1})\cong \Db(X).$$ Here,
$\Db(M,\alpha^{-1})$ is the bounded derived category of
$\Coh(M,\alpha^{-1})$ (see the next section).

Thus starting with classical untwisted K3 surfaces we naturally
end up with twisted K3 surfaces. There are other reasons to
consider twisted K3 surfaces, as has been alluded to in the
introduction, but from the point of view of moduli spaces
of sheaves on K3 surfaces this is absolutely necessary in order to
fully understand the relation  between K3 surfaces and their
moduli spaces of sheaves.

Motivated by this example, C\u{a}ld\u{a}raru formulated in
\cite{C} a conjecture that generalizes Corollary \ref{GTtrans}.
to the case of twisted K3 surfaces. In fact, the conjecture could
be verified in a number of other situations (see e.g.\
\cite{DP,Khalid}), but turned out to be wrong in general (see
\cite[Ex.\ 4.11]{HS1}).

When \cite{C} was written, generalized Calabi--Yau structures had
not been invented and the Hodge structure of twisted K3 surfaces
had not been introduced. Only the transcendental part
$T(X,\alpha)$ could be defined directly in terms of a Brauer
class $\alpha$ and Orlov's result (see Corollary \ref{GTtrans})
suggested to conjecture $\Db(X,\alpha)\cong \Db(X',\alpha')$ if
and only if $T(X,\alpha)\cong T(X',\alpha')$. However, in
contrast to the untwisted case, a Hodge iso\-metry between the
transcendental lattices of two twisted K3 surfaces does not extend
to a Hodge isometry of the full weight two Hodge structure on the
Mukai lattice. Nikulin's result does not apply any longer, as a
hyperbolic plane in the orthogonal complement does not
necessarily exist.

How the original conjecture of  C\u{a}ld\u{a}raru
has to be modified will be explained next.
\section{Twisted versions}

The question we shall deal with in this section is the following.
Suppose $X$ and $X'$ are two projective K3 surfaces endowed with
Brauer classes $\alpha$ and $\alpha'$, respectively. When does there
exist an equivalence
$$\Db(X,\alpha)\cong\Db(X',\alpha')?$$

We must confess that we are not able to deal with the question in
this generality, as an analogue of Orlov's existence result (see
Theorem \ref{OrlovExist}) in the twisted case has not yet been
proven.\footnote{In the meantime, the twisted analogue of Orlov's
result has been proved by  Canonaco and Stellari in
\cite{ConStel}. So we are talking about arbitrary exact
equivalences here.} So, whenever we deal with equivalences
between twisted derived categories in this section, we shall mean
equivalences of Fourier--Mukai type. More precisely, we will
consider
$$\xymatrix{\Phi_\kp:\Db(X,\alpha)\ar[r]^-\sim&\Db(X',\alpha')}$$
which for an object $\Db(X\times X',\alpha^{-1}\times\alpha')$ is
defined by the usual formula
$$\xymatrix{F^\bullet\ar@{|->}[r]&Rp_*(\kp\otimes^L q^*
F^\bullet).}$$ In \cite{Cal} it is explained that the usual
formalism of derived functors goes through in the twisted case.

Let us first explain the `easy' direction that has led to the
definition of the twisted Chern character in \cite{HS1} (see
Section \ref{SectTwistChern}).

\begin{proposition}
Any equivalence of Fourier--Mukai type
$$\xymatrix{\Phi_\kp:\Db(X,\alpha)\ar[r]^-\sim&\Db(X',\alpha')}$$
induces a Hodge isometry
$$\xymatrix{\Phi_\kp^H:\widetilde H
(X,\alpha,\IZ)\ar[r]^-\sim&\widetilde H(X',\alpha',\IZ).}$$
\end{proposition}

The Mukai lattice $\widetilde H(X,\alpha,\IZ)$ of a twisted K3
surface $(X,\alpha)$, has been introduced in Section
\ref{SecTwistedK3} as the Hodge structure $\widetilde
H(X,B,\IZ)$  of the generalized K3 surface given by
$\exp(B)\sigma=\sigma+B\wedge\sigma$, where $B\in H^2(X,\IQ)$
with $\exp(B^{0,2})=\alpha$. The isomorphism type of the weight
two Hodge structure $\widetilde H(X,\alpha,\IZ)$ is independent
of the choice of $B$, but for the definition of it one $B$ has to
be picked.

Thus, in order to explain the idea behind the proposition, we
need to fix B-field lifts $B\in H^2(X,\IQ)$ and $B'\in
H^2(X',\IQ)$ of $\alpha$ and $\alpha'$,  respectively. This allows us
at the same time to consider $v^{-B\oplus B'}(\kp)\in H^*(X\times
X',\IQ)$. The claimed Hodge isometry is then provided by
$$\xymatrix{\gamma\ar@{|->}[r]&p_*(v^{-B\oplus B'}(\kp).q^*\gamma).}$$

In \cite{HS1} we explain how one has to modify the arguments of
Mukai to make them work in the twisted case as well, e.g.\ why
the Mukai vector is again integral etc. The new feature in the
twisted case is that one has to make the additional choice of the
B-field lifts. In general there is no canonical lift, but if the
untwisted case is considered as a twisted case with trivial
Brauer class, then one may use the canonical lift $B=0$.

The converse of the above has been proved in \cite{HS2}. Unlike
the untwisted case, the orientation has to be incorporated in the
assertion from the beginning (see Remark \ref{Remarkderorient}).

\begin{theorem}
Suppose $(X,\alpha)$ and $(X',\alpha')$ are two projective twisted
K3 surfaces. If there exists an orientation
preserving Hodge isometry
$$\widetilde H(X,\alpha,\IZ)\cong
\widetilde H(X',\alpha',\IZ),$$ then one finds a Fourier--Mukai
equivalence
$$\xymatrix{\Phi_\kp:\Db(X,\alpha)\ar[r]^-\sim&\Db(X',\alpha').}$$
\end{theorem}

\begin{remark}\label{genericK3nosphe}
(joint work with E.\ Macr\`i and P.\ Stellari) For a generic
projective twisted K3 surface $(X,\alpha)$ there are no
$(-2)$-classes in $\widetilde H^{1,1}(X,\alpha,\IZ)$. Hence
$\Db(X,\alpha)$ does not contain any spherical objects. In this
case one can show that any Fourier--Mukai kernel
$\kp\in\Db(X\times X',\alpha^{-1}\times\alpha')$ inducing an
equivalence $\Phi_\kp:\Db(X,\alpha)\cong\Db(X',\alpha')$ is
isomorphic to a shifted sheaf $E[k]$. This is enough to conclude
that $\Phi_\kp^H$ is orientation preserving. Thus, for a generic
projective $(X,\alpha)$ the theorem reads: There exists a
Fourier--Mukai equivalence $\Db(X,\alpha)\cong\Db(X',\alpha')$ if
and only if there exists an orientation preserving Hodge isometry
$\widetilde H(X,\alpha,\IZ)\cong\widetilde H(X',\alpha',\IZ)$.
\footnote{This has now appeared in \cite{HuyMS}.}
\end{remark}

In order to prove that any orientation preserving Hodge isometry
can be lifted to a derived equivalence one tries to imitate
Orlov's proof of Theorem \ref{DGTT}. Most of the arguments go
through, but at a few crucial points the non-triviality of the
Brauer class necessitates a different approach. One is the
occasional absence of spherical objects for general twisted K3
surfaces and of the structure sheaf $\ko_X$ as an object in
$\Db(X,\alpha)$ in particular. Another one is the non-emptiness
and smoothness of certain moduli spaces of stable twisted sheaves,
which has to be assured. In particular, the following result due
to Yoshioka's plays a central r\^ole in the argument (see
\cite{Yoshioka}).

\begin{theorem}
Let $X$ be a projective K3 surface, $B\in H^2(X,\IQ)$ a B-field
and $v\in \widetilde H^{1,1}(X,B,\IZ)$ a primitive vector with
$\langle v,v\rangle=0$ and $v_0\ne0$. Then there exists a moduli
space $M(v)$ of stable (with respect to a generic polarization)
$\alpha_B$-twisted sheaves $E$ with $v^B(E)=v$ which is a
(non-empty!) K3 surface.
\end{theorem}

The existence of moduli spaces of stable twisted sheaves
has been shown in broad generality by Lieblich \cite{Lieblich}
and  Yoshioka \cite{Yoshioka}. Using the equivalence to sheaves
over Azumaya algebras,
it can also be deduced from the general results of Simpson \cite{Simpson}.
From this theorem Yoshioka deduces by standard methods the existence
of a universal $\alpha_B^{-1}\times\alpha'$-twisted sheaf $\kp$ on
$X\times M(v)$ which induces an equivalence $\Db(X,\alpha_B)\cong
\Db(M(v),\alpha')$. The latter is needed in order to imitate
iii), Example \ref{FMExK3} in the twisted case. Eventually, the
assertion is reduced to the classical Global Torelli Theorem.

\begin{remark}
The twisted version of Remark \ref{numberFM} holds true as well.
E.g.\ for a given K3 surface $X$ and a fixed Brauer class
$\alpha_0\in\Br(X)$ there are only finitely many classes
$\alpha\in\Br(X)$ such that $\Db(X,\alpha)$ is Fourier--Mukai
equivalent to $\Db(X,\alpha_0)$. See \cite[Prop.\ 3.4]{HS1}.
\end{remark}

\begin{remark}
Twisted derived equivalences have also been considered for abelian
varieties by Polishchuk (see e.g.\ \cite{KO,Polishchuk}). A
complete analogue of the untwisted results of Mukai, Orlov, and
Polishchuk has been obtained, although by methods different from
the ones in \cite{MukaiAV,OrlovAV}.
\end{remark}

\section{What's left}\label{SecEnd}
 So far we have treated `half' of the derived Global Torelli
Theorem. Staying on one K3 surface, we have up to now only tried
to determine the image of the natural representation
\begin{equation}\label{derivedequ}
\xymatrix{{\rm Aut}(\Db(X,\alpha))
\ar[r]&{\rm O}_+(\widetilde H(X,\alpha,\IZ)).}
\end{equation}
The classical Global Torelli Theorem asserts that ${\rm
Aut}(X)\to{\rm O}_+(H^2(X,\IZ))$ is injective. This is no longer
true in the derived setting, e.g.\ the shift $[2]$ is contained
in the kernel of (\ref{derivedequ}). So, the `other half' of a
derived (twisted) Global Torelli Theorem would be concerned with
the kernel of (\ref{derivedequ}). A similar question has been
asked for abelian varieties and a beautiful answer has been given
by Orlov in \cite{OrlovAV} (see \cite[Ch.\ 9]{HuyFM} for an
account of this).

For a long time the kernel of (\ref{derivedequ}) for K3 surfaces
seemed mysterious. Bridgeland's work \cite{BridgelandK3} on
stability conditions on derived categories of K3 surfaces has
changed the situation completely. We now at least have a clear
conjecture and an answer seems in reach. This would then yield
the final form of the derived, twisted Global Torelli Theorem.

Without giving any background on stability conditions, we simply
state Bridgeland's conjecture (generalized to the case of twisted
K3 surfaces).

\begin{conjecture}
For any projective twisted K3 surface $(X,\alpha)$ there exists a natural
short exact sequence
$$\xymatrix{0\ar[r]&\pi_1(\kp_0(X,\alpha))\ar[r]&
{\rm Aut}(\Db(X,\alpha))\ar[r]&{\rm O}_+(\widetilde H(X,\alpha,\IZ))
\ar[r]&1.}$$
\end{conjecture}

Before explaining what $\kp_0(X,\alpha)$ is, let us once more
recall that we actually only know that ${\rm O}_+(\widetilde
H(X,\alpha,\IZ))$, which denotes the group of all orientation
preserving Hodge isometries, is contained in the image, but we
are unable, for the time being, to show that the image is not
bigger.

Consider the set $\Delta(X,\alpha)\subset \widetilde
H^{1,1}(X,\alpha,\IZ)$ of all classes $\delta$ with
$\langle\delta,\delta\rangle=-2$. Let $\kp(X,\alpha)\subset
\widetilde H^{1,1}(X,\alpha,\IZ)\otimes_\IZ\IC$ be the open
subset of all vectors whose real and imaginary parts (in this
order) span a positive oriented plane. Then
$$\kp_0(X,\alpha):=\kp(X,\alpha)\setminus\bigcup_{\delta\in\Delta(X,\alpha)}\delta^\perp.$$

Bridgeland constructs a natural map
$\pi_1(\kp_0(X,\alpha))\to{\rm Aut}(\Db(X,\alpha))$.
(Adapting \cite{BridgelandK3} to the twisted case is
rather straightforward.)

We conclude with a few observation in the case of generic twisted
K3 surfaces $(X,\alpha)$ (joint work with E.\ Macr\`i and P.\
Stellari).

If $X$ is a generic projective K3 surface and $\alpha$ is a
generic Brauer class, then $ \Delta(X,\alpha)=\varnothing$. Thus,
$\kp_0(X,\alpha)=\kp(X,\alpha)$, whose fundamental group is $\IZ$.
This group is mapped onto the subgroup of ${\rm
Aut}(\Db(X,\alpha))$ that is spanned by the shift $[2]$. Due to
Remark \ref{genericK3nosphe} every Fourier--Mukai autoequivalence
$\Phi_\kp$ of $\Db(X,\alpha)$ has a kernel of the form $\kp\cong
E[\ell]$ for some twisted sheaf $E$ on $X\times X$ and some
$\ell\in \IZ$.

Suppose $\Phi_\kp^H={\rm id}$. Then $E_x:=E|_{\{x\}\times X}$ has
Mukai vector $(0,0,1)$. Therefore, $E_x\cong k(y)$ for some $y\in
X$ and $\ell$ must be even. As line bundle twists and
automorphisms of $X$ are all detected on cohomology, this yields

\begin{proposition}
For a generic twisted K3 surface $(X,\alpha)$ one has an exact sequence
$$\xymatrix{0\ar[r]&\IZ[2]\ar[r]&{\rm Aut}(\Db(X,\alpha))\ar[r]&{\rm O}_+(\widetilde H(X,\alpha,\IZ)))\ar[r]&1.}$$
\end{proposition}

We are not able to exclude the existence of exotic components of
the moduli space of stability conditions in the generic case,
which might be expected in general.\footnote{A proof can be found
in \cite{HuyMS}. In there, also the group of autoequivalences of
the derived category of a generic non-projective K3 surface is
determined.}



{\footnotesize
\begin{thebibliography}{mm}


\bibitem{Per} \em G\'eom\'etrie des surfaces K3: modules et p\'eriodes.
\em S\'eminaires Palaiseau. ed A.\ Beauville, J.-P.\ Bourguignon,
M.\ Demazure. Ast\'erisque 126 (1985).


\bibitem{BPV} W.\ Barth, C.\ Peters, A.\ Van de Ven
\em Compact complex surfaces. \em Springer-Verlag, Berlin (1984)

\bibitem{BO2} A.\ Bondal, D.\ Orlov
\em Semiorthogonal decomposition for algebraic varieties. \em
alg-geom/9506012.

\bibitem{BO} A.\ Bondal, D.\ Orlov
\em Reconstruction of a variety from the derived category and
groups of autoequivalences. \em Comp.\ Math.\  125  (2001),
327-344.

\bibitem{Borcea} C.\ Borcea
\em Diffeomorphisms of a K3 surface. \em Math.\ Ann.\ 275 (1986),
1-4.

\bibitem{Bridgeland} T.\ Bridgeland
\em Equivalences of Triangulated Categories and Fourier--Mukai
Transforms. \em Bull.\ London Math.\ Soc.\ 31 (1999), 25-34.


\bibitem{BM} T.\ Bridgeland, A.\ Maciocia
\em Complex surfaces with equivalent derived categories. \em
Math.\ Z.\ 236 (2001), 677--697.


\bibitem{BridgelandK3} T.\ Bridgeland
\em Stability conditions on K3 surfaces. \em math.AG/0307164.


\bibitem{BR} D.\ Burns, M.\ Rapoport
\em On the Torelli Problem for K\"ahlerian K3 Surfaces. \em Ann.\
scient.\ \'Ec.\ Norm.\ Sup.\ 8 (1975), 235-274.

\bibitem{Cal} A.\ C\u{a}ld\u{a}raru
\em Derived categories of twisted sheaves on Calabi-Yau
manifolds. \em Ph.-D. thesis Cornell (2000).

\bibitem{C} A.\ C\u{a}ld\u{a}raru
\em Non-fine moduli spaces of sheaves on $K3$ surfaces. \em
IMRN 20 (2002), 1027-1056.

\bibitem{ConStel}  A.\ Canonaco, P.\ Stellari
\em Twisted Fourier-Mukai functors. \em math.AG/0605229.


\bibitem{dJ} A.J.\ de Jong
\em A result of Gabber. \em Preprint.


\bibitem{DP} R.\ Donagi, T.\ Pantev
\em Torus fibrations, gerbes, and duality. \em to appear in:
Memoirs of the AMS. math.AG/0306213.

\bibitem{Donaldson} S.\ Donaldson
\em Polynomial invariants for smooth four-manifolds. \em Top.\ 29
(1990), 257-315.

\bibitem{FM} R.\ Friedman, J.\ Morgan
\em Smooth four-manifolds and complex surfaces. \em Erg.\ Math.\
27 (1994), Springer.


\bibitem{Gabriel} P.\ Gabriel
 \em Des cat{\'e}gories ab{\'e}liennes. \em
  Bull.\ Soc.\ Math.\ France  90  (1962), 323-448.

\bibitem{Giraud} J.\ Giraud
\em Cohomologie non-ab\'elienne. \em Springer (1971).


\bibitem{GB} A.\ Grothendieck
 \em Le groupe de Brauer II. \em In: J.~Giraud
(ed) et al.: Dix expos\'es sur la cohomologie des sch\'emas,
North-Holland, Amsterdam, (1968), 88-189.


\bibitem{Heinl} J.\ Heinloth
\em Twisted Chern Classes and $\IG_m$-gerbes. \em
C.R.\ Acad.\ Sci.\ Paris t.\ 341, 10 (2005).


\bibitem{Hit} N.\ Hitchin
\em Generalized Calabi--Yau manifolds. \em
 Q.\ J.\ Math.\  54 (2003), 281-308.

\bibitem{HLOY1} S.\ Hosono, B.H.\ Lian, K.\ Oguiso, S.-T.\ Yau
\em Autoequivalences of derived category of a K3 surface and
monodromy transformations. \em J.\ Alg.\ Geom.\ 13 (2004),
513-545.

\bibitem{Oguisonew} S.\ Hosono, B.H.\ Lian, K.\ Oguiso, S.-T.\ Yau
\em Fourier--Mukai numbers of a K3 surface. \em  CRM Proc.\ and
Lect.\ Notes 38, (2004).

\bibitem{HuyHK} D.\ Huybrechts
\em Compact hyperk\"ahler manifolds: Basic results. \em Invent.\
Math.\ 135 (1999), 63-113.

\bibitem{HuyKK} D.\ Huybrechts
\em The K\"ahler cone of a compact hyperk\"ahler manifold. \em
Math.\ Ann.\  326  (2003), 499-513.

\bibitem{HuyGen} D.\ Huybrechts
\em Generalized Calabi--Yau structures, K3 surfaces, and
B-fields. \em Int. J. Math. 16 (2005), 13-36.

\bibitem{HuyFM} D.\ Huybrechts
\em Fourier--Mukai transforms in Algebraic Geometry. \em Oxford
Mathematical Monographs (2006).

\bibitem{HL} D.\ Huybrechts, M.\ Lehn
\em The geometry of moduli spaces of shaves. \em
Aspects of Mathematics E 31, Vieweg (1997).

\bibitem{HuySch} D.\ Huybrechts, St.\ Schr\"oer
\em The Brauer group of analytic K3 surfaces. \em IMRN. 50 (2003),
2687-2698.

\bibitem{HS1} D.\ Huybrechts, P.\ Stellari
\em Equivalences of twisted K3 surfaces. \em Math.\ Ann.\ 332
(2005), 901-936.

\bibitem{HS2} D.\ Huybrechts, P.\ Stellari
\em Proof of C\u{a}ld\u{a}raru's conjecture. \em to appear in: The
13th MSJ Inter.\ Research Inst.\ Moduli Spaces and Arithmetic
Geometry, Adv.\ Stud.\ Pure Math.,  math.AG/0411541.

\bibitem{Huyab} D.\ Huybrechts
\em Derived and abelian equivalence of K3 surfaces. \em
math.AG/0604150.

\bibitem{HuyMS} D.\ Huybrechts, E.\ Macr\`i, P.\ Stellari \em Stability
conditions for generic K3 categories. \em math.AG/0608430.


\bibitem{KO} A.\ Kapustin, D.\ Orlov
\em  Vertex algebras, mirror symmetry, and D-branes: The case of
complex tori. \em  Comm.\ Math.\ Phys.\ 233 (2003), 79-136.


\bibitem{Khalid}  M.\ Khalid
\em On $K3$ correspondences. \em math.AG/0402314.

\bibitem{Lieblich} M.\ Lieblich
\em Moduli of twisted sheaves. \em math.AG/0411337.

\bibitem{Lieblich2} M.\ Lieblich
\em Sheaves and the period-index problem. \em math.AG/0511244.


\bibitem{Looij} E.\ Looijenga
 \em  A Torelli theorem for K\"ahler--Einstein $K3$ surfaces. \em
Lect.\ Notes Math.\ 894 (1981),  107-112.

\bibitem{LP} E.\ Looijenga, C.\ Peters
\em Torelli Theorems for K\"ahler K3 Surfaces. \em Comp.\ Math.\
42 (1981), 145-186.

\bibitem{Milne} J.\ Milne
\em \'Etale cohomology. \em Princeton University Press (1980).

\bibitem{Mo} D.\ Morrison
\em On K3 surfaces with large Picard number. \em Invent.\ Math.\
75, (1984) 105-121.

\bibitem{MukaiAV} S.\ Mukai
\em Duality between $D(X)$ and $D(\hat X)$ with its applications
to Picard sheaves. \em Nagoya Math.\ J.\ 81 (1981), 153-175.

\bibitem{Mu} S.\ Mukai
\em On the moduli space of bundles on K3 surfaces, I. \em In:
Vector Bundles on Algebraic Varieties, Bombay (1984).

\bibitem{Mu1} S.\ Mukai
\em Vector bundles on a $K3$ surface. \em
 Proc.\  Int.\ Congr.\ Math.\ Beijing (2002),  495-502.


\bibitem{Nikulin} V.V.\ Nikulin
\em Integral symmetric bilinear forms and some of their
applications. \em  Math USSR Izv.\ 14 (1980), 103-167.

\bibitem{Og} K.\ Oguiso
\em K3 surfaces via almost-primes. \em
 Math.\ Research Letters 9 (2002), 47-63.

\bibitem{OrlovAV} D.\ Orlov
\em Derived categories of coherent sheaves on abelian varieties
and equivalences between them. \em Izv.\ Math.\ 66 (2002),
569-594.

\bibitem{Or} D.\ Orlov
\em Equivalences of derived categories and K3 surfaces. \em J.\
Math.\ Sci.\ 84 (1997), 1361-1381.

\bibitem{Or2} D.\ Orlov
\em Derived category of coherent sheaves and equivalences between
them. \em Russ.\ Math.\ Surveys 58:3, (2003), 511-591

\bibitem{Perego} A.\ Perego
\em A Gabriel theorem for coherent twisted sheaves. \em
math.AG/0607025.

\bibitem{Ploog} D.\ Ploog
\em Groups of autoequivalences of derived categories of smooth
  projective varieties. \em PhD-thesis. FU-Berlin (2005).

\bibitem{PS} I.\ Piateckii-Shapiro, I.\ Shafarevich
\em A Torelli theorem for algebraic K3 surfaces of type K3. \em
Math.\ USSR Izv.\ 5 (1971), 547-588.

\bibitem{Pol}  A.\ Polishchuk
  \em Symplectic biextensions and generalization of the
  Fourier--Mukai transforms. \em Math.\ Res.\ Lett.\ 3 (1996),
  813-828.

\bibitem{Polishchuk} A.\ Polishschuk
\em Abelian Varieties, Theta Functions and the Fourier transform.
\em Cambridge  (2003).

\bibitem{Simpson} C.\ Simpson
\em Moduli of representations of the fundamental group of a
smooth projective variety. I.  \em  Inst.\ Hautes \'Etudes Sci.\ Publ.\
 Math.\  No.\ 79 (1994), 47-129.

\bibitem{Siu2} Y.-T.\ Siu
\em  A simple proof of the surjectivity of the period map of $K3$ surfaces. \em
Manuscripta Math.\  35  (1981), 311-321.

\bibitem{Siu} Y.-T.\ Siu
\em Every K3 surface is K\"ahler. \em
 Invent.\ math.\ 73 (1983), 139-150.

\bibitem{St} P.\  Stellari
\em Some remarks about the FM-partners of K3 surfaces with Picard
number 1 and 2. \em Geom.\ Dedicata 108 (2004), 1-13.


\bibitem{ST} P.\ Seidel, R.\ Thomas
\em Braid group actions on derived categories of coherent sheaves.
\em Duke Math.\ J.\ 108 (2001), 37-108.


\bibitem{Sz} B.\ Szendr\H{o}i
\em Diffeomorphisms and families of Fourier--Mukai transforms in
mirror symmetry. \em Applications of Alg.\ Geom.\ to Coding
Theory, Phys.\ and Comp. NATO Science Series. Kluwer (2001),
317-337.

\bibitem{Tod} A.\ Todorov
\em Applications of the K\"ahler--Einstein--Calabi--Yau metric to moduli
of $K3$ surfaces. \em
Invent.\ Math.\ 61 (1980), 251-265.

\bibitem{Verb} M.\  Verbitsky
\em Coherent sheaves on generic compact tori. \em Algebraic
structures and moduli spaces, CRM Proc.\ Lecture Notes, 38,
(2004), 229-247.


\bibitem{W} A.\ Weil
 \em Final report on contract AF 18(603)-57. \em
Collected works (1979).

\bibitem{Yoshioka} K.\ Yoshioka
\em Moduli spaces of twisted sheaves on projective varieties. \em
to appear in: The 13th MSJ Inter.\ Research Inst.\ Moduli Spaces
and Arithmetic Geometry, Adv.\ Stud.\ Pure Math., AMS,
math.AG/0411538.


\end{thebibliography}}
\end{document}